\documentclass[11pt]{article}
\usepackage[francais]{babel}
\usepackage[latin1]{inputenc}
\usepackage[T1]{fontenc}
\usepackage{a4,amsmath,amsfonts}
\title{\'Equidistribution des sous-variétés de petite hauteur}
\author{Pascal Autissier}

\begin{document}

\maketitle

\newcommand{\D}{\displaystyle}

{\bf Abstract}\\

In this paper, the equidistribution theorem of Szpiro-Ullmo-Zhang about
sequences of small points in an abelian variety is extended to the case of
sequences of higher dimensional subvarieties. A quantitative version of this
result is also given.\\

{\bf\Large Introduction}\\

Dans ce texte, on appelle {\bf variété} sur un corps $K$ tout schéma intègre
et projectif sur $K$.

Commençons par rappeler les travaux de Szpiro, Ullmo et Zhang sur
l'équidistribution des petits points ({\it cf} \cite{SUZ}):\\

{\bf Définition:} Soit $V$ une variété sur $\mathbb{Q}$. Une suite
$(Y_n)_{n\ge0}$ de fermés intègres de $V$ est dite {\bf générique} dans $V$
lorsque pour tout fermé $Z\neq V$, l'ensemble
$\{n\in\mathbb{N}\ |\ Y_n\subset Z\}$ est fini.\\

Soit $K$ un corps de nombres de degré $N$. Notons $G_K$ l'ensemble des
plongements $\sigma:K\hookrightarrow\mathbb{C}$. Soit $A$ une variété abélienne
sur $K$ de dimension $d$. Remarquons que $A(\mathbb{C})$ est alors l'union
disjointe des groupes $A_\sigma(\mathbb{C})$. Soit $\mu$ la mesure de
probabilité sur $A(\mathbb{C})$ égale sur chaque $A_\sigma(\mathbb{C})$ à la
mesure de Haar de masse $1/N$.\\

Soit $\hat{\rm h}_L$ la hauteur de Néron-Tate relativement à un
faisceau inversible $L$ symétrique et ample sur $A$. Soit $(P_n)_{n\ge0}$ une
suite de points fermés de $A$. Pour tout $n\ge0$, on pose
$k_n=\deg(P_n)=[k(P_n):\mathbb{Q}]$.\\

Lorsque $Z$ est un fermé réduit de $A_\mathbb{C}$ purement de codimension $q$,
on désigne par $\delta_Z$ le courant d'intégration sur $Z(\mathbb{C})$
(c'est un $(q;q)$-courant sur $A(\mathbb{C})$).\\

Szpiro, Ullmo et Zhang ont démontré le théorème d'équirépartition suivant:\\

{\bf Théorème 0:} {\it Supposons que la suite $(P_n)_{n\ge0}$ est générique
dans $A$ et que $\D\lim_{n\rightarrow+\infty}\hat{\rm h}_L(P_n)=0$. Alors la
suite de mesures $\D\Bigl(\frac{1}{k_n}\delta_{P_{n\mathbb{C}}}\Bigr)_{n\ge0}$
converge faiblement vers $\mu$.}\\

Ce résultat a permis à Ullmo \cite{Ull} et Zhang \cite{Zh3} de prouver la
conjecture de Bogomolov ({\it cf} aussi \cite{Abb}).\\

On étend ici le théorème 0 aux suites génériques de sous-variétés:\\

Soit $\|\ \|$ une ``métrique du cube'' ({\it cf} \cite{Mor}) sur
$L_\mathbb{C}$. Notons $\omega$ la $(1;1)$-forme de courbure de
$(L_\mathbb{C};\|\ \|)$ sur $A(\mathbb{C})$ et posons $k=\deg_L(A)$.
Remarquons que l'on a $\omega^{\wedge d}=k\mu$ sur $A(\mathbb{C})$.\\

Soit $(Y_n)_{n\ge0}$ une suite de fermés intègres de $A$ de dimension $p$. On
pose $k_n=\deg_L(Y_n)$ pour tout $n\ge0$. On prouve le résultat suivant:\\

{\bf Théorème (4.1):} {\it Supposons que la suite $(Y_n)_{n\ge0}$ est générique
dans $A$ et que $\D\lim_{n\rightarrow+\infty}\hat{\rm h}_L(Y_n)=0$. Alors la
suite de mesures
$\D\Bigl(\frac{1}{k_n}\omega^{\wedge p}\wedge\delta_{Y_{n\mathbb{C}}}\Bigr)_{n\ge0}$ converge faiblement vers $\mu$.}\\

La démonstration s'inspire de celle du théorème 0. Elle se place donc dans le
cadre de la théorie d'Arakelov ({\it cf} \cite{BGS}), et utilise le théorème de
``Hilbert-Samuel arithmétique'' d\^u à Gillet et Soulé \cite{GiSo}
({\it cf} aussi Abbes et Bouche \cite{AbBo}).\\

Remarquons en particulier que si $(Y_n)_{n\ge0}$ est générique dans $A$ et si
$U$ est un ouvert non vide de $A(\mathbb{C})$ disjoint de $Y_{n\mathbb{C}}$
pour tout $n$ assez grand, alors on a
$\D\ \liminf_{n\rightarrow+\infty}\hat{\rm h}_L(Y_n)>0$. On peut donner une
version quantitative de ce résultat:\\

{\bf Théorème (6.1):} {\it On suppose la suite $(Y_n)_{n\ge0}$ générique dans
$A$. Soit $U_\eta$ une boule ouverte de $A(\mathbb{C})$ de rayon
$\eta\in]0;\eta_0[$ telle que $U_\eta$ soit disjoint de son conjugué complexe.
On suppose que $U_\eta$ et $Y_{n\mathbb{C}}$ sont disjoints pour tout $n\ge0$.
Alors on a
$\D\ \liminf_{n\rightarrow+\infty}\hat{\rm h}_L(Y_n)\ge c_0\eta^{2d+2}$
($\eta_0$ et $c_0$ sont des constantes ne dépendant que de $(K;A;L;p)$).}\\

Par ailleurs, Zhang a proposé une généralisation de la conjecture de
Bogomolov ({\it cf} \cite{Zh2}):\\

Soient $V$ une variété lisse sur $\mathbb{Q}$ et $L$ un faisceau inversible
ample sur $V$. Supposons que l'on a un morphisme fini $f:V\rightarrow V$ et un
isomorphisme $L^{\otimes m}\simeq f^*L$, où $m$ est un entier $>1$.\\

{\bf Définition:} Un fermé intègre $Y$ de $V$ est dit {\bf prépériodique}
lorsque l'ensemble $\{f^n(Y)\ ;\ n\in\mathbb{N}\}$ est fini.\\

Lorsque l'on a de telles données, Zhang a montré comment obtenir une hauteur
canonique $\hat{\rm h}_L$ généralisant la construction de Néron-Tate. En
particulier, la fonction $\hat{\rm h}_L$ est positive; si un fermé intègre $Y$
de $V$ est prépériodique, alors on a $\hat{\rm h}_L(Y)=0$; et un point fermé
$P$ est prépériodique si et seulement si $\hat{\rm h}_L(P)=0$.\\

Dans cette situation, la conjecture de Zhang s'énonce ainsi:\\

{\bf Conjecture (i):} Soit $Y$ un fermé intègre de $V$. Si
$\hat{\rm h}_L(Y)=0$, alors $Y$ est prépériodique.\\

Zhang a montré que $\hat{\rm h}_L(Y)=0$ si et seulement si $Y$ contient une
suite de points fermés $(P_n)_{n\ge0}$ générique dans $Y$ telle que
$\D\lim_{n\rightarrow+\infty}\hat{\rm h}_L(P_n)=0$. En particulier, si $Y$
contient un ensemble Zariski-dense de points fermés prépériodiques, alors $Y$
est conjecturalement prépériodique.

Ce dernier énoncé, qui ne fait plus intervenir de hauteur, rappelle une
conjecture de dynamique complexe:\\

{\bf Conjecture (ii):} Soient $M$ une variété lisse sur $\mathbb{C}$ et
$f:M\rightarrow M$ un morphisme fini. On suppose qu'il existe un faisceau
inversible $L$ ample sur $M$ tel que $L^{\otimes m}\simeq f^*L$
avec $m>1$. Si un fermé intègre $Y$ de $M$ contient un ensemble Zariski-dense
de points prépériodiques, alors $Y$ est prépériodique.\\

Cette conjecture ``s'applique'' en particulier au cas où $M$ est l'espace
projectif $\mathbb{P}^d_\mathbb{C}$ et $f$ un morphisme fini de degré $>1$
(on prend alors $L={\cal O}(1)$).

Par ailleurs, lorsque $M$ est une variété abélienne sur $\mathbb{C}$ et $f$ la
multiplication par un entier $>1$, la conjecture (ii) est un théorème de
Raynaud ({\it cf} \cite{Ray}).\\

Pour attaquer la conjecture (i), on pourrait commencer par prouver un résultat
d'équidistribution des petits points analogue au théorème 0. La mesure limite
$\mu$ serait alors la mesure à l'équilibre du système dynamique
$(V_\mathbb{C};f_\mathbb{C})$.

La difficulté de cette approche est que l'on est amené, semble-t-il, à utiliser
un théorème de Hilbert-Samuel arithmétique sans hypothèse de positivité de la
courbure, qui n'est connu que pour les surfaces arithmétiques ({\it cf}
propositions 3.3.3 et 4.1.4 de \cite{Aut}).\\

Je remercie Romain Dujardin et Charles Favre pour des discussions intéressantes
concernant la conjecture (ii). Je remercie également Antoine Chambert-Loir pour
ses suggestions concernant la section 5. Je remercie Serge Cantat pour de
fructueuses conversations et pour ses conseils concernant ce papier. Enfin, je
remercie Emmanuel Ullmo pour l'inspiration qu'il m'a procurée.\\

\section{Définitions et notations}

Soient $M$ une variété lisse sur $\mathbb{C}$ et $L$ un faisceau inversible
sur $M$.\\

{\bf Définition:} Une {\bf métrique} sur $L$ est une famille $\|\ \|$ de normes
hermitiennes sur $L$ variant de manière continue sur $M(\mathbb{C})$. Posons
$\widehat{L}=(L;\|\ \|)$.

On note $\omega_{\widehat{L}}$ le $(1;1)$-courant de courbure de
$\widehat{L}$ sur $M(\mathbb{C})$; il est caractérisé par la propriété
suivante: pour toute section rationnelle $s$ non nulle de $L$, on a l'égalité
$\D\ \omega_{\widehat{L}}=\delta_{{\rm div}(s)}-\frac{i}{\pi}\partial\overline{\partial}\ln\|s\|$.\\

{\bf Définition:} Soit $\|\ \|$ une métrique sur $L$; posons
$\widehat{L}=(L;\|\ \|)$. On dit que $\|\ \|$ est {\bf p.s.h.} lorsque le
courant $\omega_{\widehat{L}}$ est positif, autrement dit lorsque pour tout
ouvert $U$ de $M$ et tout $s\in\Gamma(U;L)$ ne s'annulant pas sur $U$, la
fonction $\ -\ln\|s\|$ est plurisousharmonique sur $U(\mathbb{C})$.\\

{\bf Définitions:} Une {\bf variété arithmétique} est un schéma $X$ intègre,
projectif et plat sur $\mathbb{Z}$, tel que la fibre générique $X_\mathbb{Q}$
soit lisse sur $\mathbb{Q}$. Remarquons que $X_\mathbb{C}$ est alors l'union
disjointe de variétés lisses sur $\mathbb{C}$.

Un {\bf faisceau inversible hermitien} sur $X$ est un couple
$\widehat{\cal L}=({\cal L};\|\ \|)$, formé d'un faisceau inversible ${\cal L}$
sur $X$ et d'une métrique $\|\ \|$ sur ${\cal L}_\mathbb{C}$ invariante par
conjugaison complexe.\\

{\bf Remarque:} Soit $Y$ un schéma intègre et projectif sur $\mathbb{Z}$. On a
deux possibilités:\\
- $Y$ est plat et surjectif sur $B_0={\rm Spec}(\mathbb{Z})$; $Y$ est alors
dit horizontal;\\
- $Y$ est au-dessus d'un point fermé $b=p\mathbb{Z}$ de $B_0$ ({\it ie} $Y$
est une variété sur $\mathbb{F}_p$); $Y$ est alors dit vertical.\\

\section{Théorie d'Arakelov}

Soient $X$ une variété arithmétique de dimension (absolue) $d$, et
$\widehat{\cal L}=({\cal L};\|\ \|)$ un faisceau inversible hermitien sur $X$
dont la métrique $\|\ \|$ est ${\rm C}^\infty$. Pour $0\le p\le d$, on désigne
par ${\rm Z}_p(X)$ le groupe des combinaisons $\mathbb{Z}$-linéaires de fermés
intègres (de $X$) de dimension $p$.\\

Faisons quelques rappels de la théorie des hauteurs développée par Bost, Gillet
et Soulé:\\

Pour $0\le p\le d$, on définit ({\it cf} \cite{GiSo} p. 485, \cite{BGS} p. 933)
le groupe de Chow arithmétique de dimension $p$, noté $\widehat{\rm CH}_p(X)$.
Le degré arithmétique définit une application $\mathbb{Z}$-linéaire
$\ \widehat{\deg}:\widehat{\rm CH}_0(X)\rightarrow\mathbb{R}$. Pour
$1\le p\le d$, on construit la première classe de Chern arithmétique
$\ \widehat{\rm c}_1(\widehat{\cal L}):\widehat{\rm CH}_p(X)\rightarrow\widehat{\rm CH}_{p-1}(X)$, qui est $\mathbb{Z}$-linéaire.\\

Soient maintenant $p\in\{0;\cdots;d\}$ et $D\in{\rm Z}_p(X)$. Soit $g$ un
courant de Green pour $D_\mathbb{C}$. Le réel
$${\rm h}_{\widehat{\cal L}}(D)=\widehat{\deg}\Bigl[\widehat{\rm c}_1(\widehat{\cal L})^p(D;g)-(0;\omega_{\widehat{\cal L}_\mathbb{C}}^{\wedge p}\wedge g)\Bigr]$$
ne dépend pas du choix de $g$. On l'appelle la {\bf hauteur d'Arakelov} de $D$
relativement à $\widehat{\cal L}$.\\

Soit $Y$ un fermé intègre horizontal de dimension $p$ tel que
$\deg_{{\cal L}_\mathbb{Q}}(Y_\mathbb{Q})>0$. La {\bf hauteur
normalisée} de $Y$ relativement à $\widehat{\cal L}$ est le réel
$${\rm h}'_{\widehat{\cal L}}(Y)=\frac{{\rm h}_{\widehat{\cal L}}(Y)}{\deg_{{\cal L}_\mathbb{Q}}(Y_\mathbb{Q})p}\quad.$$

{\bf Proposition 2.1:} {\it Soient $Y$ un fermé intègre de $X$ de dimension
$p\ge1$, $n$ un entier $\ge1$, et $s$ une section rationnelle non nulle de
${\cal L}_{|Y}^{\otimes n}$ (sur $Y$). Alors on a l'égalité
$${\rm h}_{\widehat{\cal L}}(Y)n={\rm h}_{\widehat{\cal L}}({\rm div}(s))-\int_{Y(\mathbb{C})}\ln\|s_\mathbb{C}\|\omega_{\widehat{\cal L}_\mathbb{C}}^{p-1}$$
(On convient que l'intégrale est nulle si $Y$ est vertical).}\\

{\it Démonstration:} C'est la proposition 3.2.1 (iv) de \cite{BGS} p. 949.
$\square$\\

Le résultat suivant est un analogue arithmétique du théorème de
Hilbert-Samuel:\\

{\bf Théorème 2.2:} {\it Supposons ${\cal L}$ ample et la métrique $\|\ \|$
p.s.h.. Soit $\epsilon>0$. Alors pour tout entier $n$ assez grand, il existe
une section globale $s$ non nulle de ${\cal L}^{\otimes n}$ telle que
$\D\ \max_{X(\mathbb{C})}\ln\|s_\mathbb{C}\|\le\epsilon n-{\rm h}'_{\widehat{\cal L}}(X)n$.}\\

{\it Démonstration:} {\it cf} théorème 9 de \cite{GiSo} p. 539 ou corollaire du
théorème principal de \cite{AbBo}. $\square$\\

\section{Hauteurs}

\subsection{Hauteurs de Weil}

Soient $V$ une variété lisse sur $\mathbb{Q}$ et $L$ un faisceau inversible
ample sur $V$. Notons $V^\diamond$ l'ensemble des fermés intègres de $V$.\\

{\bf Définition:} Un {\bf modèle entier} de $(V;L)$ est un couple
$(X;\widehat{\cal L})$, formé d'une variété arithmétique $X$ et d'un faisceau
inversible hermitien $\widehat{\cal L}=({\cal L};\|\ \|)$ sur $X$ tels que
$(X_\mathbb{Q};{\cal L}_\mathbb{Q})$ soit isomorphe à $(V;L)$, que ${\cal L}$
soit ample sur $X$, que la métrique $\|\ \|$ soit ${\rm C}^\infty$, et que la
courbure $\omega_{\widehat{\cal L}_\mathbb{C}}$ soit définie positive sur
$X(\mathbb{C})$.

Remarquons que pour tout entier $e$ assez grand, le couple $(V;L^{\otimes e})$
admet un modèle entier.\\

La définition suivante étend aux sous-variétés la notion classique de hauteur
de Weil sur les points:\\

{\bf Définition:} Une application $\hat{h}:V^\diamond\rightarrow\mathbb{R}$
est une {\bf hauteur de Weil} relativement à $L$ lorsqu'il existe un entier
$e\ge1$, un modèle entier $(X;\widehat{\cal L})$ de $(V;L^{\otimes e})$ et un
réel $C$ tels que
$$\forall Y\in V^\diamond\quad\Bigl|\hat{h}(Y)-\frac{1}{e}{\rm h}'_{\widehat{\cal L}}(\overline{Y})\Bigr|\le C$$
($\overline{Y}$ désigne l'adhérence de $Y$ dans $X$).\\

La proposition suivante montre entre autres que cette définition ``ne dépend
pas du choix de $(e;X;\widehat{\cal L})$'':\\

{\bf Proposition 3.1:} {\it Soit $\hat{h}$ une hauteur de Weil relativement à
$L$. Alors pour tout entier $e\ge1$ et tout modèle entier
$(X;\widehat{\cal L})$ de $(V;L^{\otimes e})$, il existe un réel $C$ tel que
$\quad\D\forall Y\in V^\diamond\quad\Bigl|\hat{h}(Y)-\frac{1}{e}{\rm h}'_{\widehat{\cal L}}(\overline{Y})\Bigr|\le C$.

En outre, la fonction $\hat{h}$ est minorée sur $V^\diamond$.}\\

{\it Démonstration:} La première partie de l'énoncé est une reformulation de la
proposition 3.2.2 de \cite{BGS} p. 950. La deuxième partie se déduit de la
remarque (iii) de \cite{BGS} p. 954. $\square$\\

\subsection{Hauteur canonique}

Soient $V$ une variété lisse sur $\mathbb{Q}$ et $L$ un faisceau inversible
ample sur $V$. On suppose que l'on a un morphisme fini $f:V\rightarrow V$ et un
isomorphisme $\alpha:L^{\otimes m}\xrightarrow{\sim}f^*L$, où $m$ est un
entier $>1$.\\

\`A partir de ces données, Zhang construit une hauteur de Weil particulière
relativement à $L$ par un procédé dynamique:\\

{\bf Proposition 3.2:} {\it Il existe une unique hauteur de Weil
$\hat{\rm h}_L$ relativement à $L$ telle que pour tout $Y\in V^\diamond$, on
ait l'égalité $\ \hat{\rm h}_L(f(Y))=\hat{\rm h}_L(Y)m$.

De plus, on a les propriétés suivantes:\\
- Si $\hat{h}$ est une hauteur de Weil relativement à $L$, alors la suite de
fonctions $\D\Bigl(\frac{1}{m^n}\hat{h}\circ f^n\Bigr)_{n\ge0}$ converge
uniformément sur $V^\diamond$ vers la fonction $\hat{\rm h}_L$;\\
- La fonction $\hat{\rm h}_L$ est positive sur $V^\diamond$;\\
- Si $Y$ est un fermé intègre prépériodique de $V$, alors on a
$\hat{\rm h}_L(Y)=0$;\\
- Un point fermé $P$ de $V$ est prépériodique si et seulement si
$\hat{\rm h}_L(P)=0$.}\\

{\it Démonstration:} C'est une reformulation du théorème 2.4 de \cite{Zh2}
p. 292. $\square$\\

L'application $\hat{\rm h}_L$ est appelé la {\bf hauteur canonique}
relativement à $(L;f)$. Remarquons qu'elle ne dépend pas du choix de
l'isomorphisme $\alpha$.\\

\section{\'Equidistribution}

Soit $K$ un corps de nombres de degré $N$. Soient $A$ une variété abélienne sur
$K$ de dimension $d$ et $0_A\in A(K)$ sa section neutre. Pour tout entier non
nul $c$, on note $[c]:A\rightarrow A$ le morphisme de multiplication par $c$.\\

Soit $L$ un faisceau inversible symétrique ({\it ie} $[-1]^*L\simeq L$) et
ample sur $A$. On fixe un isomorphisme $0_A^*L\simeq{\cal O}_B$ sur
$B={\rm Spec}(K)$. Par le théorème du cube, on en déduit naturellement un
isomorphisme $[c]^*L\simeq L^{\otimes c^2}$ sur $A$ pour chaque entier
$c\ge2$.\\

La hauteur canonique $\hat{\rm h}_L$ relativement à $(L;[c])$ ne dépend pas
du choix de l'entier $c\ge2$; c'est la {\bf hauteur de Néron-Tate} relativement
à $L$ ({\it cf} aussi \cite{Phi} et \cite{Gub}).\\

On fixe un entier $c\ge2$. D'après \cite{Mor} p. 50-52, il existe une unique
métrique $\|\ \|$ de classe ${\rm C}^\infty$ sur $L_\mathbb{C}$ telle que
l'isomorphisme $[c]^*L\simeq L^{\otimes c^2}$ devienne une isométrie ({\it ie}
une ``métrique du cube'').

Notons $\omega$ la courbure de $(L_\mathbb{C};\|\ \|)$; posons $k=\deg_L(A)$
et $\D\mu=\frac{1}{k}\omega^d$. Pour tout plongement
$\sigma:K\hookrightarrow\mathbb{C}$, la mesure $\mu_\sigma$ est alors la mesure
de Haar de masse $1/N$ sur $A_\sigma(\mathbb{C})$.\\

Soit $(Y_n)_{n\ge0}$ une suite de fermés intègres de $A$ de dimension $p$.
On pose $k_n=\deg_L(Y_n)$ pour tout $n\ge0$.\\

{\bf Théorème 4.1:} {\it Supposons que la suite $(Y_n)_{n\ge0}$ est générique
dans $A$ et que $\D\lim_{n\rightarrow+\infty}\hat{\rm h}_L(Y_n)=0$. Alors la
suite de mesures
$\D\Bigl(\frac{1}{k_n}\omega^p\delta_{Y_{n\mathbb{C}}}\Bigr)_{n\ge0}$ converge
faiblement vers $\mu$.}\\

{\it Démonstration:} Posons
$\D\nu_n=\frac{1}{k_n}\omega^p\delta_{Y_{n\mathbb{C}}}$ pour tout entier
$n\ge0$. Il s'agit de montrer que pour toute fonction $\phi$ continue sur
$A(\mathbb{C})$, la suite numérique
$\D\Bigl(\int_{A(\mathbb{C})}\phi\nu_n\Bigr)_{n\ge0}$ converge vers
$\D\int_{A(\mathbb{C})}\phi\mu$.

Il suffit en fait de le vérifier pour $\phi$ invariante par conjugaison
complexe et ${\rm C}^\infty$ sur $A(\mathbb{C})$. Il existe alors un réel
$t_0>0$ tel que pour tout $t\in[-t_0;t_0]$, la $(1;1)$-forme
$\D\omega_t=\omega+\frac{it}{\pi}\partial\overline{\partial}\phi$ soit positive
sur $A(\mathbb{C})$. Soient $t\in[-t_0;t_0]$ et $\epsilon>0$.\\

On choisit un entier $e\ge1$ et un modèle
$(X;\widehat{\cal L})$ de $(A;L^{\otimes e})$ tels que la métrique sur
${\cal L}$ soit la métrique du cube $\|\ \|^{\otimes e}$.

Pour tout $Y\in A^\diamond$, on pose
$\D\hat{h}(Y)=\frac{1}{e}{\rm h}'_{\widehat{\cal L}}(\overline{Y})$, où
$\overline{Y}$ désigne l'adhérence de $Y$ dans $X$. L'application $\hat{h}$ est
par définition une hauteur de Weil relativement à $L$. D'après la proposition
3.2, il existe $n_1\ge0$ tel que
$\D\Bigl|\hat{h}\circ[c^{n_1}]-c^{2n_1}\hat{\rm h}_L\Bigr|\le c^{2n_1}\epsilon$
sur $A^\diamond$. On pose $e_1=c^{2n_1}e$.\\

On construit un modèle $(X';\widehat{\cal L}')$ de $(A;L^{\otimes e_1})$ de la
manière suivante:\\

On désigne par $f:X'\rightarrow X$ la normalisation de $X$ par le morphisme
$[c^{n_1}]:A\rightarrow A$, de sorte que l'on a $f_\mathbb{Q}=[c^{n_1}]$ sur
$X'_\mathbb{Q}\simeq A$. On pose alors
$\widehat{\cal L}'=({\cal L}';\|\ \|')=f^*\widehat{\cal L}$.\\

Par la formule de projection ({\it cf} proposition 3.2.1 (iii) de \cite{BGS}
p. 949), on a, pour tout $Y\in A^\diamond$, la relation
$\hat{h}([c^{n_1}](Y))e={\rm h}'_{\widehat{\cal L}'}(\overline{Y}')$, où
$\overline{Y}'$ désigne l'adhérence de $Y$ dans $X'$. On en déduit:
$$\forall Y\in A^\diamond\quad\Bigl|\hat{\rm h}_L(Y)-\frac{1}{e_1}{\rm h}'_{\widehat{\cal L}'}(\overline{Y}')\Bigr|\le\epsilon\quad.\qquad(1)$$

On pose $\|\ \|'_t=\|\ \|'\exp(-e_1t\phi)$ et
$\widehat{\cal L}'_t=({\cal L}';\|\ \|'_t)$. Par multilinéarité, on a la
formule 
$${\rm h}_{\widehat{\cal L}'_t}(X')-{\rm h}_{\widehat{\cal L}'}(X')=(d+1)e_1^{d+1}kt^2Q(t)+(d+1)e_1^{d+1}t\int_{A(\mathbb{C})}\phi\omega^d\quad,$$
où $\D Q(T)=\sum_{j=1}^d\frac{{\rm C}_{d+1}^{j+1}}{(d+1)k}T^{j-1}\int_{A(\mathbb{C})}\phi\Bigl(\frac{i}{\pi}\partial\overline{\partial}\phi\Bigr)^j\omega^{d-j}$ est un polyn\^ome en $T$.\\

Sachant que $\hat{\rm h}_L(A)=0$, on en déduit à l'aide de (1) la minoration
$$\frac{1}{e_1}{\rm h}'_{\widehat{\cal L}'_t}(X')\ge t^2Q(t)-\epsilon+t\int_{A(\mathbb{C})}\phi\mu\quad.\qquad(2)$$

Maintenant, appliquons le théorème 2.2 (Hilbert-Samuel arithmétique): il
existe $n_2\ge1$ et $s\in\Gamma(X';{\cal L}'^{\otimes n_2})-\{0\}$ tels que
$\D\ \max_{A(\mathbb{C})}\ln\|s_\mathbb{C}\|'_t\le\epsilon n_2-{\rm h}'_{\widehat{\cal L}'_t}(X')n_2$.

Posons $Z={\rm div}(s_\mathbb{Q})$. Puisque la suite $(Y_n)_{n\ge0}$ est
générique dans $A$, il existe un entier $n_3\ge0$ tel que
$\ \forall n\ge n_3\ Y_n\not\subset Z$. Alors, d'après la proposition 2.1, on a
(pour tout $n\ge n_3$):
$$\begin{array}{rcl}
\D{\rm h}_{\widehat{\cal L}'}(\overline{Y_n}')n_2-{\rm h}_{\widehat{\cal L}'}({\rm div}(s_{|\overline{Y_n}'}))&=&\D-e_1^p\int_{Y_n(\mathbb{C})}(n_2e_1t\phi+\ln\|s_\mathbb{C}\|'_t)\omega^p\\
&\ge&\D[{\rm h}'_{\widehat{\cal L}'_t}(X')-\epsilon]e_1^pn_2k_n-e_1^{p+1}n_2k_nt\int_{A(\mathbb{C})}\phi\nu_n\ .\\
\end{array}$$

Majorons le premier membre: En utilisant (1) et la positivité de
$\hat{\rm h}_L$, on obtient d'une part l'inégalité ${\rm h}_{\widehat{\cal L}'}({\rm div}(s_{|\overline{Y_n}'}))\ge-e_1^{p+1}pn_2k_n\epsilon$ et d'autre part
l'inégalité ${\rm h}_{\widehat{\cal L}'}(\overline{Y_n}')\le[\hat{\rm h}_L(Y_n)+\epsilon](p+1)e_1^{p+1}k_n$. On en déduit à l'aide de la
minoration (2) que l'on a (pour tout $n\ge n_3$):
$$(p+1)\hat{\rm h}_L(Y_n)+t\int_{A(\mathbb{C})}\phi\nu_n\ge t^2Q(t)-(2p+3)\epsilon+t\int_{A(\mathbb{C})}\phi\mu\quad.\qquad(3)$$

On fait tendre $n$ vers $+\infty$ puis $\epsilon$ vers 0 dans ce qui précède;
sachant que $\D\lim_{n\rightarrow+\infty}\hat{\rm h}_L(Y_n)=0$, on trouve
l'inégalité
$$\liminf_{n\rightarrow+\infty}\Bigl(t\int_{A(\mathbb{C})}\phi\nu_n\Bigr)\ge t^2Q(t)+t\int_{A(\mathbb{C})}\phi\mu\quad.$$

Finalement, en faisant tendre $t$ vers 0:\\
- Par valeurs supérieures, on obtient $\D\ \liminf_{n\rightarrow+\infty}\Bigl(\int_{A(\mathbb{C})}\phi\nu_n\Bigr)\ge\int_{A(\mathbb{C})}\phi\mu$;\\
- Par valeurs inférieures, on obtient $\D\ \limsup_{n\rightarrow+\infty}\Bigl(\int_{A(\mathbb{C})}\phi\nu_n\Bigr)\le\int_{A(\mathbb{C})}\phi\mu$.\\

On en déduit le résultat. $\square$\\

\section{Variante arakelovienne}

Soient $X$ une variété arithmétique de dimension $d$ et
$\widehat{\cal L}=({\cal L};\|\ \|)$ un faisceau inversible hermitien sur $X$.
On suppose que ${\cal L}$ est ample sur $X$, que la métrique $\|\ \|$ est
${\rm C}^\infty$, et que la courbure
$\omega=\omega_{\widehat{\cal L}_\mathbb{C}}$ est définie positive sur
$X(\mathbb{C})$. Le couple $(X;\widehat{\cal L})$ est donc un modèle entier de
$(X_\mathbb{Q};{\cal L}_\mathbb{Q})$.\\

Posons $k=\deg_{{\cal L}_\mathbb{Q}}(X_\mathbb{Q})$ et
$\D\mu=\frac{1}{k}\omega^{d-1}$. Soit $(Y_n)_{n\ge0}$ une suite de fermés
intègres de $X_\mathbb{Q}$ de dimension $p$. On pose
$k_n=\deg_{{\cal L}_\mathbb{Q}}(Y_n)$ pour tout $n\ge0$.\\

On fait l'hypothèse $(*)$ suivante: pour tout fermé intègre $Y$ de
$X_\mathbb{Q}$ de dimension $p-1$, on a
${\rm h}'_{\widehat{\cal L}}(\overline{Y})\ge{\rm h}'_{\widehat{\cal L}}(X)$.

Remarquons que cette hypothèse est automatiquement vérifiée lorsque $p=0$.\\

{\bf Proposition 5.1:} {\it Supposons que la suite $(Y_n)_{n\ge0}$ est
générique dans $X_\mathbb{Q}$ et que
$\D\lim_{n\rightarrow+\infty}{\rm h}'_{\widehat{\cal L}}(\overline{Y_n})={\rm h}'_{\widehat{\cal L}}(X)$. Alors la suite de mesures
$\D\Bigl(\frac{1}{k_n}\omega^p\delta_{Y_{n\mathbb{C}}}\Bigr)_{n\ge0}$ converge
faiblement vers $\mu$.}\\

{\it Démonstration:} Posons
$\D\nu_n=\frac{1}{k_n}\omega^p\delta_{Y_{n\mathbb{C}}}$ pour tout entier
$n\ge0$. Soit $\phi$ une fonction invariante par conjugaison complexe et
${\rm C}^\infty$ sur $X(\mathbb{C})$. Il existe un réel $t_0>0$ tel que pour
tout $t\in[-t_0;t_0]$, la $(1;1)$-forme $\D\omega_t=\omega+\frac{it}{\pi}\partial\overline{\partial}\phi$ soit positive sur $X(\mathbb{C})$.\\

Soient $t\in[-t_0;t_0]$ et $\epsilon>0$. On pose $\|\ \|_t=\|\ \|\exp(-t\phi)$
et $\widehat{\cal L}_t=({\cal L};\|\ \|_t)$. Par multilinéarité, on a la
relation
$${\rm h}_{\widehat{\cal L}_t}(X)-{\rm h}_{\widehat{\cal L}}(X)=dkt^2Q(t)+dt\int_{X(\mathbb{C})}\phi\omega^{d-1}\quad,$$
où $Q$ est une fonction polyn\^omiale.\\

D'après le théorème 2.2, il existe $n_1\ge1$ et $s\in\Gamma(X;{\cal L}^{\otimes n_1})-\{0\}$ tels que $\D\max_{X(\mathbb{C})}\ln\|s_\mathbb{C}\|_t\le\epsilon n_1-{\rm h}'_{\widehat{\cal L}_t}(X)n_1$.

Posons $Z={\rm div}(s_\mathbb{Q})$. La suite $(Y_n)_{n\ge0}$ étant générique
dans $X_\mathbb{Q}$, il existe un entier $n_2\ge0$ tel que
$\ \forall n\ge n_2\ Y_n\not\subset Z$. Alors, en appliquant la
proposition 2.1, on obtient (pour tout $n\ge n_2$):
$$\begin{array}{rcl}
\D{\rm h}_{\widehat{\cal L}}(\overline{Y_n})n_1-{\rm h}_{\widehat{\cal L}}({\rm div}(s_{|\overline{Y_n}}))&=&\D-\int_{Y_n(\mathbb{C})}(n_1t\phi+\ln\|s_\mathbb{C}\|_t)\omega^p\\
&\ge&\D[{\rm h}'_{\widehat{\cal L}_t}(X)-\epsilon]n_1k_n-n_1k_nt\int_{X(\mathbb{C})}\phi\nu_n\quad.\\
\end{array}$$

Gr\^ace à l'hypothèse $(*)$, on a la minoration
$\ {\rm h}_{\widehat{\cal L}}({\rm div}(s_{|\overline{Y_n}}))\ge pn_1k_n{\rm h}'_{\widehat{\cal L}}(X)$.
On en déduit que pour tout $n\ge n_2$, on a
$$(p+1){\rm h}'_{\widehat{\cal L}}(\overline{Y_n})+t\int_{X(\mathbb{C})}\phi\nu_n\ge(p+1){\rm h}'_{\widehat{\cal L}}(X)+t^2Q(t)-\epsilon+t\int_{X(\mathbb{C})}\phi\mu\quad.$$

On conclut comme dans la démonstration du théorème 4.1. $\square$\\

\section{Version quantitative}

Soient $K$ un corps de nombres, $A$ une variété abélienne sur $K$ de dimension
$d$, et $L$ un faisceau inversible symétrique et ample sur $A$. Soit $\|\ \|$
une métrique du cube sur $L_\mathbb{C}$. On note $\omega$ la courbure de
$(L_\mathbb{C};\|\ \|)$ et on pose $k=\deg_L(A)$.\\

On fixe un plongement $\sigma:K\hookrightarrow\mathbb{C}$ et un
isomorphisme $A_\sigma(\mathbb{C})\simeq\mathbb{C}^d/\Gamma$ de groupes
analytiques (où $\Gamma$ est un réseau de $\mathbb{C}^d$) tel que la
$(1;1)$-forme $\omega_\sigma$ s'écrive
$\D\omega_\sigma=\frac{i}{2}\sum_{j=1}^d{\rm d}z_j\wedge{\rm d}\overline{z_j}$
dans $\mathbb{C}^d/\Gamma$.

On munit $\mathbb{C}^d$ de la norme hermitienne canonique, notée
$|\!|\ |\!|$. Lorsque $z\in\mathbb{C}^d$, on note ici $\dot{z}$ l'image de $z$
dans $\mathbb{C}^d/\Gamma$. On munit $\mathbb{C}^d/\Gamma$ de la distance $D$
induite par $|\!|\ |\!|$, {\it ie}
$\D D(\dot{x};\dot{y})=\min_{\gamma\in\Gamma}|\!|x-y+\gamma|\!|$. Par ailleurs,
on pose $\D\eta_0=\frac{1}{2}\min_{\gamma\in\Gamma-\{0\}}|\!|\gamma|\!|$.\\

Soit $(Y_n)_{n\ge0}$ une suite de fermés intègres de $A$ de dimension $p$.
Posons $\D c_0=\frac{2\pi^d}{(p+1)(d+1)k}$. Soit $U_\eta$ une boule ouverte de
$A_\sigma(\mathbb{C})$ de rayon $\eta\in]0;\eta_0[$. Si $\sigma$ est réel, on
suppose $U_\eta$ disjoint de son conjugué complexe.\\

{\bf Théorème 6.1:} {\it On suppose que la suite $(Y_n)_{n\ge0}$ est générique
dans $A$ et que $Y_{n\mathbb{C}}$ est disjoint de $U_\eta$ pour tout $n\ge0$.
Alors on a
$\D\ \liminf_{n\rightarrow+\infty}\hat{\rm h}_L(Y_n)\ge c_0\eta^{2d+2}$.}\\

{\it Démonstration:} Soient $\dot{z_0}$ le centre de $U_\eta$ et $U'_\eta$ le
conjugué complexe de $U_\eta$. Soit $\psi:[0;\eta_0[\rightarrow\mathbb{R}$ une
fonction ${\rm C}^\infty$ vérifiant:\\
- La fonction $\psi$ est nulle sur $[\eta;\eta_0[$;\\
- On a $\psi'\ge-1$ et $\psi''\ge0$ sur $[0;\eta_0[$.\\

Soit $\phi:A(\mathbb{C})\rightarrow\mathbb{R}$ la fonction caractérisée par les
propriétés suivantes:\\
- Si $|\!|z-z_0|\!|<\eta$, alors
$\D\phi(\dot{z})=\frac{\pi}{2}\psi\Bigl(|\!|z-z_0|\!|^2\Bigr)$;\\
- La fonction $\phi$ est invariante par conjugaison complexe;\\
- La fonction $\phi$ est nulle en dehors de $U_\eta\cup U'_\eta$.\\

La $(1;1)$-forme $\D\omega_1=\omega+\frac{i}{\pi}\partial\overline{\partial}\phi$ est alors positive sur $A(\mathbb{C})$. On reprend avec $t_0=t=1$ la
démonstration du théorème 4.1, jusqu'à l'inégalité (3):
$$\forall n\ge n_3\quad(p+1)\hat{\rm h}_L(Y_n)+\int_{A(\mathbb{C})}\phi\nu_n\ge Q_1-(2p+3)\epsilon\quad,$$
où $\D Q_1=\frac{1}{k}\sum_{j=0}^d\int_{A(\mathbb{C})}\phi\omega_1^j\omega^{d-j}$.\\

L'intégrale dans le premier membre est nulle puisque $Y_{n\mathbb{C}}$ et
$U_\eta$ sont disjoints. En faisant tendre $n$ vers $+\infty$ puis $\epsilon$
vers 0 dans l'inégalité précédente, on obtient donc la minoration
$\D\ \liminf_{n\rightarrow+\infty}\hat{\rm h}_L(Y_n)\ge\frac{Q_1}{p+1}$.\\

Par ailleurs, un calcul que j'épargne au lecteur montre que
$$Q_1=\frac{4\pi^d}{k}\int_0^\eta\Bigl[1-\Bigl(1+\psi'(r^2)\Bigr)^{d+1}\Bigr]r^{2d+1}{\rm d}r\quad.$$

On a donc
$$\liminf_{n\rightarrow+\infty}\hat{\rm h}_L(Y_n)\ge c_0\int_0^\eta\Bigl[1-\Bigl(1+\psi'(r^2)\Bigr)^{d+1}\Bigr](2d+2)r^{2d+1}{\rm d}r\quad.$$

Maintenant, soit $\psi_0:[0;\eta_0[\rightarrow\mathbb{R}$ la fonction définie
par $\psi_0(\rho)=\eta-\rho$ si $\rho<\eta$ et $\psi_0(\rho)=0$ si
$\rho\ge\eta$. En faisant tendre $\psi$ vers $\psi_0$ convenablement, on en
déduit le résultat. $\square$\\

\section{Cas des courbes}

Soient $V$ une courbe lisse sur $\mathbb{Q}$, $L$ un faisceau inversible
ample sur $V$, $f:V\rightarrow V$ un morphisme fini, et
$\alpha:L^{\otimes m}\xrightarrow{\sim}f^*L$ un isomorphisme, où $m$ est un
entier $>1$.\\

Notons $\mu$ la mesure à l'équilibre du système dynamique
$(V_\mathbb{C};f_\mathbb{C})$. Soit $(P_n)_{n\ge0}$ une suite de points fermés
de $V$. Pour tout $n\ge0$, on pose $k_n=[k(P_n):\mathbb{Q}]$.\\

Citons pour mémoire la variante suivante de la proposition 4.1.4 de
\cite{Aut}:\\

{\bf Proposition 7.1:} {\it Supposons que la suite $(P_n)_{n\ge0}$ est
générique dans $V$ et que $\D\lim_{n\rightarrow+\infty}\hat{\rm h}_L(P_n)=0$.
Alors la suite de mesures
$\D\Bigl(\frac{1}{k_n}\delta_{P_{n\mathbb{C}}}\Bigr)_{n\ge0}$ converge
faiblement vers $\mu$.}\\

\ \\

{\small Pascal Autissier. $\qquad$ pascal.autissier@univ-rennes1.fr\\
I.R.M.A.R., Université de Rennes I, campus de Beaulieu, 35042 Rennes cedex,
France.}

\end{document}